\documentclass{birkjour}
 \newtheorem{thm}{Theorem}[section]
 
 \newtheorem{lem}[thm]{Lemma}
 \newtheorem{prop}[thm]{Proposition}
 \theoremstyle{definition}
 \newtheorem{defn}[thm]{Definition}
 \theoremstyle{remark}
 \newtheorem{rem}[thm]{Remark}
 
 \numberwithin{equation}{section}
 \usepackage{mathabx}
\usepackage{amssymb}
\usepackage{mathrsfs}
\usepackage{hyperref}
\usepackage{enumerate}
\usepackage{physics}
\usepackage[notcite, final, notref]{showkeys}
\usepackage{dsfont}

\newcommand{\e}{\varepsilon}

\newcommand{\w}{\omega}

\usepackage{tikz}
\usetikzlibrary{arrows,shapes,positioning,shadows,trees,matrix}
\usepackage{xcolor}
\title[]{Analytic continuation of time in Brownian motion. Stochastic distributions approach}

\author[L.D. Abreu]{Luis Daniel Abreu}
\address{NuHAG University of Vienna\\
  Vienna\\
  Austria}
\email{abreuluisdaniel@gmail.com}

\author[D. Alpay]{Daniel Alpay}
\address{Schmid College of Science and Technology \\
Chapman University\\
One University Drive
Orange, California 92866\\
USA}
\email{alpay@chapman.edu}

\author[T. Georgiou]{Tryphon Georgiou}
\address{Department of Mechanical and Aerospace Engineering\\
  University of californila at Irvine,\\
  Irvine\\
CA, USA} 
\email{tryphon@uci.edu}

\author[P. Jorgensen]{Palle Jorgensen}
\address{Department of Mathematics, 14 MLH\\The University
of Iowa, Iowa City, Iowa 52242-1419\\ USA}
\email{palle-jorgensen@uiowa.edu}

\thanks{      \section*{Acknowledgments}
           Daniel Alpay thanks the Foster G. and Mary McGaw Professorship in Mathematical Sciences, which supported this research}

\begin{document}
         
\begin{abstract}
  With the use of Hida’s white noise space theory and spaces of stochastic distributions, we present a detailed analytic continuation theory for classes of Gaussian processes, with focus here on Brownian motion. For the latter, we prove and make use a priori bounds, in the complex plane,
  for the Hermite functions; as well as a new approach to stochastic distributions. This in turn allows us to present (in Section \ref{white-z}) an explicit formula for an analytically continued white noise process, realized this way in complex domain. With the use of the Wick product, we then apply our complex
  white noise analysis in Section \ref{sec-7} in a derivation of a new realization of Hilbert space-valued stochastic integrals.
\end{abstract}


\maketitle

\section{Introduction}
\setcounter{equation}{0}

Using Hida's white noise space theory and an associated Gelfand triple of stochastic test functions and stochastic distributions, we construct a
Gaussian process $(B_z)$ indexed by the complex numbers, and which extends in a natural way the Brownian
motion $(B_t)$ indexed by $\mathbb R$. In this section we explain the general strategy. The various notions are recalled in the
paper itself.\smallskip

Consider Hida's white noise space $(\mathscr S^{\prime}_{\mathbb R},\mathcal C,P)$, where $\mathscr S^{\prime}_{\mathbb R}$ denote the
  space of real-valued tempered distributions and $\mathcal C$ is the associated cylinder sigma-algebra. See Section \ref{sec-w-n}.
  Furthermore, let $\zeta_0,\zeta_1,\ldots$ denote the normalized Hermite functions (see \eqref{hnnorm}), and let $Z_0,Z_1,\ldots$ be the
  family of independent $N(0,1)$
  variables in $\mathbf L_2(\mathscr S^{\prime}_{\mathbb R},\mathcal C,P)$ defined by \eqref{zn}. The Brownian motion can be constructed as
    \begin{equation}
B_t(\w)=\sum_{n=0}^\infty\left( \int_0^t \zeta_n(u)du\right) Z_n(\w),\quad t\in\mathbb R,
    \end{equation}
where the convergence is in $\mathbf L_2(\mathscr S^{\prime}_{\mathbb R},\mathcal C,P)$; see e.g. \cite{new_sde}. This suggests to define for $z\in \mathbb C$

\begin{equation}
  \label{BMz}
B_z(\w)=\sum_{n=0}^\infty\left( \int_{[0,z]} \zeta_n(u)du\right) Z_n(\w).
    \end{equation}
    The values of the integrals $\int_{[0,z]} \zeta_n(u)du$ do not depend on the path linking $0$ to $z$ since the Hermite functions are entire.
    On the other hand, Mehler's formula gives for $|\epsilon|<1$
    \begin{equation}
      \label{mehler678}
      \begin{split}
      \sum_{n=0}^\infty \e^n \left(\int_{[0,z]}\zeta_n(u)du\right)\overline{\left(\int_{[0,w]}\zeta_n(v)dv\right)}&=\\
      &\hspace{-3.cm}   =   \frac{1}{\sqrt{2\pi}}\frac{1}{\sqrt{1-\e^2}}\iint_{[0,z]\times[0,w]}e^{-\frac{(1+\e^2)(u^2+\overline{v}^2)-4u\overline{v}\e}
        {2(1-\e^2)} }     dudv,
      \end{split}
    \end{equation}
    and will not converge, in particular for $z=w=iT$, $T>0$, as $\e\longrightarrow 1$; see Section \ref{section-8}.
    Thus, for $z\not\in\mathbb R$, the sum  \eqref{BMz} is {\sl a priori} formal and does not seem to converge in
    $\mathbf L_2(\mathscr S^{\prime}_{\mathbb R},\mathcal C,P)$. It converges in a space of
    stochastic distributions, denoted here $\mathfrak F_{dist}(\mathbf 2)$ (see \eqref{Fa} and Definition \ref{2}) with corresponding space of stochastic test functions
    $\mathfrak F_{test}(\mathbf 2)$ (see \eqref{Fa-1} and Definition \ref{2}).\\

The spaces $\mathfrak F_{test}(a)$  and $\mathfrak F_{dist}(a)$  are
    of the kind introduced in a special case by V\r{a}ge \cite{vage96b,vage96} and Kondratiev \cite{new_sde,MR1407325} and
    studied by one of the authors with Guy Salomon in a series of papers which includes \cite{MR3029153,MR3404695}. It follows that the covariance of
    the Brownian motion can be defined as an operator from $\mathfrak F_{test}(a)$ into $\mathfrak F_{dist}(a)$ in the following way:
    \begin{equation}
(B_z\otimes \overline{B_w})f=B_z\langle \overline{B_w},f\rangle
      \end{equation}
      where the brackets denote the duality between the space of stochastic distributions and the space of stochastic test functions.
      We refer to \cite{lifs,MR521027,MR626346} for other instances where covariance operators are defined as positive operators
      from a topological vector space into its anti-dual.\smallskip
      
    Applications of the Kondratiev
    space of stochastic distributions include also \cite{MR2610579} for the theory of linear systems and models for stochastic processes,
    and \cite{MR3231624} for the non-commutative white noise.\smallskip
    
    The spaces of stochastic distributions considered here are nuclear, and are inductive limits of increasing sequences of Hilbert spaces with
    decreasing norms; they are topological algebras where a sequence of inequalities relating the norms of these Hilbert spaces hold. Using the theory of
    nuclear spaces as exposed in \cite{GS2_english,MR0435834} and a bound on the absolute value of Hermite polynomials
    given in \cite{MR1041203} (see \eqref{bound-1} below) these inequalities allow to define stochastic integrals associated to $B_z$. These spaces are not metrizable, but the underlying structure allows to work with sequences, as explain in the sequel.\\

    One of the early motivations for studying complex time-Brownian motion derives in part from quantum physics: As explained in, for example, the papers by Ed Nelson \cite{MR161189,MR214150}, the corresponding analytic continuation issues, in their original form, are linked to variants of the Feynman-Kac question.
   While the existing literature deals with other issues related to Brownian motion in
    the complex domain (see e.g.
    \cite{MR3456934,MR4666291,MR3882989}), it will be clear from the comments above, and details in our paper below, that our present results go in a different direction. Nonetheless, it is our hope that our results will shine new light on the existing literature.\smallskip

    The paper consists of ten sections of which this introduction is the first.
Before getting to our main results, we must first present the necessary framework: white noise spaces, new stochastic test-function spaces and their duals and bounds on Hermite functions. This is done in Sections \ref{sec-2}-\ref{sec-sto}. More precisely, bounds for complex values of the Hermite functions are developped in Section \ref{sec-2}. We also
    recall in that section a few facts on Fr\'echet nuclear spaces and their dual. The main aspects of Hida's white noise space theory are reviewed in Section \ref{sec-w-n}, while stochastic test functions and stochastic distributions are discussed in Seection \ref{sec-sto}. Then, in Section \ref{sec-5}
    we include our results for complex Brownian motion (introduced in sections 1 and 2); and in Section \ref{white-z}, our associated complex white noise construction,
    Theorems \ref{analyticity-thm}-\ref{6.4}. With this we are then able to present a new class of stochastic integrals in our complex framework,
    Theorem \ref{7.1}. The last two sections deal with a discussion of real vs complex time via Mehler’s formula, and a more general family of complex stochastic processes,     Theorem \ref{9.2}, and a possible strategy to obtain It\^o's formula in the present setting.
    
    \section{Preliminaries}
    \label{sec-2}

We review some facts and definitions on Hermite functions and on dual of nuclear Fr\'echet spaces; see e.g. \cite{MR0205028,MR37:726} for more information on the topic.\\ 

{\bf Hermite polynomials and Hermite functions.}
We need results involving Hermite polynomials taken form \cite{MR1041203} and \cite{new_sde}.
The Hermite polynomials can be defined in two different ways (physicists versus probabilists), that is
\begin{equation}
  H_n(z)= (-1)^ne^{z^2}\left(\frac{\partial}{\partial z}\right)^ne^{-z^2},\quad n=0,1,2,\ldots
\end{equation}
(see e.g. \cite[(3) p. 431]{MR1502747}, \cite[p. 302]{sansone}, \cite[(5.5.3) p. 106]{szego}), used in particular in \cite{MR1041203}, versus
\begin{equation}
  \label{zhang-123}
  h_n(z)= (-1)^ne^{\frac{z^2}{2}}\left(\frac{\partial}{\partial z}\right)^ne^{-\frac{z^2}{2}},\quad n=0,1,2,\ldots
\end{equation}
used in \cite{new_sde}. The two definitions are related by
\begin{equation}
  H_n(x)=2^{n/2}h_n(\sqrt{2}x),\quad n=0,1,\ldots
\end{equation}

We denote by $\zeta_0,\zeta_1,\ldots$ the normalized Hermite functions
(see e.g. \cite[p. 4]{thangavelu1993lectures}),
\begin{equation}
\label{hnnorm}
\zeta_n(z)=\frac{1}{\pi^{1/4}2^{n/2}(n!)^{1/2}}(-1)^{n}e^{-\frac{z^2}{2}}H_n(z),\quad n=0,1,\ldots
\end{equation}
They form an orthonormal basis of the space $\mathbf L_2(\mathbb R,dx)$.
They satisfy in particular the following relationships (see e.g. \cite[(A.17) p. 255]{MR1851117} for \eqref{recur} and \cite[(A.26) p. 312]{Hida_BM} for \eqref{recur-1})
\begin{eqnarray}
  \label{recur}
  \zeta_n^\prime(x)&=&-x\zeta_n(x)+\sqrt{2n}\zeta_{n-1}(x),\quad n=1,2,\ldots\\
  \zeta_n^\prime(x)&=&x\zeta_n(x)-\sqrt{2(n+1)}\zeta_{n+1}(x),\quad n=1,2,\ldots\\
  x\zeta_n(x)&=&\sqrt{\frac{n}{2}}\zeta_{n-1}(x)+\sqrt{\frac{n+1}{2}}\zeta_{n+1}(x).
                 \label{recur-1}
  \end{eqnarray}

We will need bounds for $H_n$ in the complex plane in order to obtain related bounds for the Hermite functions. Eijndhoven and Meyers assert in \cite[(1.2)]{MR1041203} that
  \begin{equation}
    \label{bound-1}
    |H_n(z)|\le 2^{n/2}(n!)^{1/2}e^{\sqrt{2n}|z|}
  \end{equation}
  where it is crucial to note for the sequel that $z$ is complex, and not confined to the real line. No proof appears in \cite{MR1041203}.
  A proof can be found at the URL \cite{4839872},  where a paper of Rusev \cite{MR1808423} is also mentioned. We will use Rusev's paper, who gives the following bound
  \begin{equation}
    |H_n(z)|\le
        K_n e^{x^2}\cosh((2n+1)^{\frac{1}{2}}y)
  \end{equation}
  where
  \begin{equation}
    \label{Kn}
K_n=\sqrt[4]{\frac{2e}{\pi}}(\Gamma(2n+1))^{\frac{1}{4}}(2n+1)^{-\frac{n}{2}-\frac{1}{4}}e^{\frac{n}{2}},
    \end{equation}
    on Hermite polynomials, to give a bound on the Hermite functions. Before proceeding with these bounds we make a remark.
    In \cite[(5) p. 431]{MR1502747} one can find the following bound:
\begin{prop}
  Let $|z|\le R$. Then
  \begin{equation}
    \label{hille-345}
|H_n(z)|\le n!e^{R|z|+\frac{n}{2}}\frac{1}{R^n}.
    \end{equation}
    \end{prop}
    This bound could be used in the present setting, but is much rougher than the one obtained from Rusev's paper \cite{MR1808423}; see Remark \ref{rem-321-123}.

 \begin{lem}
    The following bound holds in the complex plane for the normalized Hermite functions:
  \begin{equation}
    \label{ineq-k}
|\zeta_n(z)|\le    \sqrt[4]{\frac{2e^3}{\pi^2}}e^{\frac{|z|^2}{2}+x^2+|y|}e^{2n|y|},\quad n=1,2,\ldots
              \end{equation}
      \end{lem}
\begin{proof}
    We recall (see \cite[p. 54]{MR0228020}) that
    \begin{equation}
      \label{feller56}
    \sqrt{2\pi}n^{n+\frac{1}{2}}e^{-n}e^{\frac{1}{12n+1}}\le n!\le     \sqrt{2\pi}n^{n+\frac{1}{2}}e^{-n}e^{\frac{1}{12n}},\quad n=1,2,\ldots
    \end{equation}
    and hence
    \begin{equation}
      \label{labek}
      \sqrt{2\pi}n^{n+\frac{1}{2}}e^{-n}\le n!\le     \sqrt{2\pi}n^{n+\frac{1}{2}}e^{-n}e,\quad n=1,2,\ldots
      \end{equation}

 In the expression \eqref{Kn} for $K_n$ we remark that $\Gamma(2n+1)=(2n)!$ and that $2n+1\le 2n$. Hence, taking into account \eqref{labek}
      we have
      \[
      \begin{split}
        K_n&\le\sqrt[4]{\frac{2e}{\pi}}\left(\sqrt{2\pi}(2n)^{2n+\frac{1}{2}}e^{-2n}e\right)^{1/2}\frac{1}{(2n)^{\frac{1}{4}}(2n)^{\frac{n}{2}}}
          e^{\frac{n}{2}}\\
        &= \sqrt[4]{\frac{2e}{\pi}}\frac{2^{\frac{1}{4}}\pi^{\frac{1}{4}}2^{n+\frac{1}{4}}n^{n+\frac{1}{4}}e^{-n}e^{\frac{1}{2}}
          e^{\frac{n}{2}}}{2^{\frac{1}{4}}n^{\frac{1}{4}}2^{\frac{n}{2}}n^{\frac{n}{2}}}\\
        &=\sqrt[4]{4e^3}2^{\frac{n}{2}}n^{\frac{n}{2}}e^{\frac{-n}{2}}.
      \end{split}
      \]
Thus,
      \[
      \begin{split}
        |\zeta_n(z)|&\le      \frac{1}{\pi^{\frac{1}{4}}2^{\frac{n}{2}}(n!)^{-\frac{n}{2}}}e^{\frac{|z|^2}{2}}
        \sqrt[4]{4e^3}2^{\frac{n}{2}}n^{\frac{n}{2}}e^{\frac{n}{2}}e^{x^2}\cosh((2n+1)^{\frac{1}{2}}y)\\
        &\le        \frac{1}{\pi^{\frac{1}{4}}2^{\frac{n}{2}}
          \left((2\pi)^{\frac{1}{2}}  n^{n+\frac{1}{2}} e^{-n}\right)     ^{\frac{1}{2}}}e^{\frac{|z|^2}{2}}
        \sqrt[4]{4e^3}2^{\frac{n}{2}}n^{\frac{n}{2}}e^{-\frac{n}{2}}e^{x^2}\cosh((2n+1)^{\frac{1}{2}}y)\\
        &\le \frac{\sqrt[4]{4e^3}}{\pi^{\frac{1}{2}}2^{\frac{1}{4}}}e^{\frac{|z|^2}{2}+x^2}
        e^{(2n+1)|y|}\\
        &=\sqrt[4]{\frac{2e^3}{\pi^2}}e^{\frac{|z|^2}{2}+x^2+|y|}
        e^{2n|y|},
      \end{split}
      \]
      where we have used in particular the inequalities
      \[
n^{\frac{1}{2}}\ge 1\quad{\rm and}\quad \cosh((2n+1)^{\frac{1}{2}}y)\le e^{(2n+1)^{\frac{1}{2}}|y|}\le e^{(2n+1)|y|}
      \]
      to go from the second-to-last line to the last one.
  \end{proof}

  Similarly we have:

    \begin{lem}
    The following bounds hold in the complex plane for the first derivatives of  the Hermite functions:
  \begin{eqnarray}
    \label{ineq-k-2}
    |\zeta_n^\prime(z)|&\le&
                             \sqrt[4]{\frac{2e^3}{\pi^2}}e^{\frac{|z|^2}{2}+x^2+|y|}
                             \left(1+\frac{\sqrt{2n}}{e^{2|y|}}\right)e^{2n|y|}\\
        \label{ineq-k-3}
    &\le&
     \sqrt[4]{\frac{2e^3}{\pi^2}}e^{\frac{|z|^2}{2}+x^2+|y|}\left(1+\sqrt{2n}\right)e^{2n|y|}.
      \end{eqnarray}
      \end{lem}
    \begin{proof}

  Equation \eqref{recur} extends analytically to the complex plane and
  it follows that
  \begin{equation}
    \begin{split}
      |\zeta_n^\prime(z)|&\le |z|\cdot|\zeta_n(z)|+\sqrt{2n}|\zeta_{n+1}(z)|\\
      &\le \sqrt[4]{\frac{2e^3}{\pi^2}}e^{\frac{|z|^2}{2}+x^2+|y|}\left(e^{2n|y|}+\sqrt{2n}e^{(2n-2)|y|}\right)\\
      &      \le \sqrt[4]{\frac{2e^3}{\pi^2}}e^{\frac{|z|^2}{2}+x^2+|y|}\left(e^{2n|y|}+\sqrt{2n}e^{(2n-2)|y|+1}\right)\\
&\le  \sqrt[4]{\frac{2e^3}{\pi^2}}e^{\frac{|z|^2}{2}+x^2+|y|}\left(1+\frac{\sqrt{2n}}{e^{2|y|}}\right)e^{2n|y|}
        \end{split}
      \end{equation}
where we have used \eqref{ineq-k} twice.
\end{proof}

\begin{defn}
  We set
  \begin{equation}
    c_z=\sqrt[4]{\frac{2e^3}{\pi^2}}e^{\frac{|z|^2}{2}+x^2+|y|},
    \label{K_z}
    \end{equation}
    and, for $R>0$,
    \begin{equation}
      \label{C_R}
      C_R=\max_{|z|\le R}c_z.
      \end{equation}
  \end{defn}
  
  In the sequel, to define the topological algebras which define the spaces of stochastic distributions,  we will
  need superexponential sequences of weights, i.e. sequences  $(a_n)$ of positive numbers indexed by $\mathbb N_0$ and which satisfy
  \begin{equation}
\label{anaman+m}
    a_na_m\le a_{n+m},\quad n,m=0,1,\ldots
  \end{equation}
  See Theorem \ref{dallas}.\\

  \begin{rem}
    \label{rem-321-123}
    Going back to the bound \eqref{hille-345} we obtain by definition of $\zeta_n$ (see \eqref{hnnorm})
    \[
      \begin{split}
        |\zeta_n(z)|&=\frac{1}{\pi^{1/4}2^{n/2}(n!)^{1/2}}|e^{-\frac{z^2}{2}}|\cdot |H_n(z)|\\
        &\le \frac{e^{\frac{R^2}{2}}}{\pi^{1/4}2^{n/2}(n!)^{1/2}}(n!)e^{R|z|+\frac{n}{2}}\frac{1}{R^n}\quad ({\rm using\,\, \eqref{hille-345}})\\
        &\le \frac{e^{\frac{R^2}{2}+R|z|}}{\pi^{1/4}}
        \left(\frac{e^{\frac{1}{2}}}{2^{\frac{1}{2}}R}\right)^n(n!)^{1/2}.
      \end{split}
    \]
    Let
    \begin{equation}
      \label{ann}
      a_n=\left(\frac{e^{\frac{1}{2}}}{2^{\frac{1}{2}}R}\right)^n(n!)^{1/2}.
    \end{equation}
    Since
    \[
     {n+m \choose n}=\frac{(n+m)!}{n!m!}\ge 1
    \]
    we have
    \[
      n!m!\le(n+m)!\quad {\rm and \,\,\ so},\,\       \sqrt{n!}\sqrt{m!}\le\sqrt{(n+m)!}
    \]
    and $(a_n)$ defined by \eqref{ann} satisfies \eqref{anaman+m}. So some of the following arguments work also for this sequence, but the bound is much bigger, leading to a bigger space where the arguments hold.
    \end{rem}
  
  {\bf Dual of a nuclear Fr\'echet space:}

  We consider the projective limit of a sequence $\mathcal H_p$ of decreasing Hilbert spaces with increasing norms, and assume that it is nuclear.
  We denote $\mathcal H_{p}^\prime=\mathcal H_{-p}$.
  The dual $\bigcup_{p=1}^\infty \mathcal H_{-p}$ endowed with the inductive topology is not metrizable (and so the topology is not determined by sequences), but
  it has the following properties, which are the key to the arguments in this paper:

  \begin{thm}
    \label{gelfand-1}
    In the above setting and notation:
  \begin{enumerate}
  \item The strong and and the inductive topologies coincide.

\item A set is bounded in $\bigcup_{p=1}^\infty \mathcal H_{-p}$ in the topology of
    $\bigcup_{p=1}^\infty \mathcal H_{-p}$ if and only it is bounded in one of the spaces $\mathcal H_{-p}$, with respect to the topology of
    $\mathcal H_{-p}$.

\item A set is compact in $\bigcup_{p=1}^\infty \mathcal H_{-p}$ in the topology of
    $\bigcup_{p=1}^\infty \mathcal H_{-p}$ if and only it is compact in one of the spaces $\mathcal H_{-p}$, with respect to the topology of
    $\mathcal H_{-p}$.

  \item
    A sequence of elements in $\bigcup_{p=1}^\infty \mathcal H_{-p}$ converges in the inductive topology if and only if it
    lies in one of the $\mathcal H_{-p}$ and converges there in the $\mathcal H_{-p}$ topology. 

      \end{enumerate}
    \end{thm}

    The last property is important when one studies continuity of functions from a metric space into $\bigcup_{p=1}^\infty \mathcal H_{-p}$.

    \section{The white noise space}
    \setcounter{equation}{0}
    \label{sec-w-n}

    As underlying probability space we will chose Hida' white noise space $(\mathscr S^{\prime}_{\mathbb R},\mathcal C,P)$   where (as above)
    $\mathscr S^{\prime}_{\mathbb R}$ denotes the
    space of real valued tempered distributions (the dual of the space $\mathscr S_{\mathbb R}$ of real-valued Schwartz functions)    and $\mathcal C$ the associated cylinder sigma-algebra. 
  The reason is that the associated $\mathbf L_2(\mathscr S^{\prime}_{\mathbb R},\mathcal C,P)$  can be identified in a natural way with the Fock space
  associated to the space $\ell_2(\mathbb N)$, i.e the reproducing kernel Hilbert space with reproducing kernel
  \[
K(\mathbf z,\mathbf w)=e^{\langle \mathbf z,\mathbf w\rangle_{\ell_2(\mathbb N)}}.
    \]
    See e.g. \cite{new_sde}.\\

    We briefly review some aspects of Hida's white noise space needed here, and in particular
ecall the construction of the measure $P$,  and the main aspects of Hida's theory relevant here. See \cite{MR562914,Hi90,MR1387829} for further informaton\smallskip
    
    The first step in Hida's theory is to apply the Bochner-Minlos theorem to the function
\begin{equation}
    \varphi(s)=e^{-\frac{\|s\|_2^2}{2}}
\end{equation}
where $s$ varies in  $\mathscr S_{\mathbb R}$, and $\|s\|_2$ denotes the norm in the Lebesgue space
$\mathbf L_2(\mathbb R,dx)$. Since $\mathscr S_{\mathbb R}$ is nuclear, the Bochner-Minlos theorem insures that there exists a positive measure
on the cylinder sigma-sigma algebra on the dual $\mathscr S^{\prime}_{\mathbb R}$ (the real-valued tempered distributions in one variable) such that
\begin{equation}
  \label{iso-12345}
e^{-\frac{\|s\|_2^2}{2}}=\int_{\mathscr S^{\prime}_{\mathbb R}}e^{i\langle s,\w\rangle} dP(\w),
\end{equation}
where the brackets denote the duality between $\mathscr S_{\mathbb R}$ and $\mathscr S^{\prime}_{\mathbb R}$, and more generally
between $\mathscr S$ and $\mathscr S^{\prime}$.\\

It follows from \eqref{iso-12345} that the random variable
\begin{equation}
  Q_s(\w)=\langle s,\w\rangle
\end{equation}
is a centered Gaussian with covariance $\|s\|^2_2$. We denote
\begin{equation}  
  \label{zn}
  Z_n=Q_{\zeta_n},\quad n=1,2,\ldots,
\end{equation}
and remark that the map $s\mapsto Q_s$ extends to an isometry from $\mathbf L_2(\mathbb R,dx)$ into
$\mathbf L_2(\mathscr S_{\mathbb R}^\prime,\mathcal C,P)$.\smallskip

Next we introduce an orthogonal basis of $\mathbf L_2(\mathscr S^{\prime}_{\mathbb R},\mathcal C,P)$, and to conclude, the Wick product. The latter gives the
motivation to introduce the stochastic distributions in the present setting.

\begin{defn}
We denote by $\ell$ the set of infinite sequences
\[
  (\alpha_1,\alpha_2,\ldots)
  \]
  indexed  by $\mathbb N$, and with values in $\mathbb N_0$, and for which $\alpha_j\not=0$ for at most a finite number of indices.
\label{ell}
\end{defn}
  \begin{prop} (see e.g. \cite[Theorem 2.2.3 p. 24]{new_sde}) With $h_0,h_1,h_2,\ldots$ denote the Hermite polynomials as defined by \eqref{zhang-123}, and $\zeta_1,\zeta_2,\ldots$ denote the
    Hermite functions the random variables
    \begin{equation}
      H_\alpha(\w)=\prod_{n=1}^\infty h_{\alpha_n}(Q_{\zeta_n})
    \end{equation}
    form an orthonormal basis of $\mathbf L_2(\mathscr S^{\prime}_{\mathbb R},\mathcal C,P)$, and
    \begin{equation}
      \langle H_\alpha, H_\beta\rangle_P=\delta_{\alpha,\beta}\alpha!.
    \end{equation}
    \label{prop-123}
  \end{prop}

  The Wick product is defined by
  \begin{equation}
    H_\alpha\star H_\beta=H_{\alpha+\beta},\quad \alpha,\beta\in\ell.
  \end{equation}
  It is thus the convolution, or Cauchy product, on $\ell$. The Wick product is not a law of composition on
  $\mathbf L_2(\mathscr S^{\prime}_{\mathbb R},\mathcal C,P)$, and we imbed the latter in a topological algebra of a special type (called strong algebra in
  \cite{MR3404695}) and in which the Wick product is stable.\\

\section{Stochastic test functions and stochastic distributions}
\setcounter{equation}{0}
\label{sec-sto}
This section follows in particular \cite{vage1}.
Consider a sequence of numbers $(a_n)_{n\in\mathbb N_0}$ satisfying $a_0=1$ and
    \begin{equation}
\label{super1}
      a_n a_m\le a_{n+m}
\end{equation}
(i.e. a superexponential sequence with $a_0=1$) and we assume that there exists $d\in\mathbb N$ such that
\begin{equation}
  \sum_{n=1}^\infty a_n^{-d}<\infty.
\label{super}
\end{equation}

For $\alpha\in\ell$  (see Definition \ref{ell}), we denote

\begin{equation}
b_\alpha=\prod a_n^{\alpha_n}.
      \end{equation}
From \eqref{super} we have
    \begin{equation}
b_\alpha b_\beta\le b_{\alpha+\beta},\quad \alpha,\beta\in\ell.
\end{equation}

\begin{lem}
  In the above notations,
  \begin{equation}
\sum_{\alpha\in\ell}b_\alpha^{-d}<\infty.
    \end{equation}
  \end{lem}

  \begin{proof}
    \[
  \sum_{\alpha\in\ell}b_\alpha^{-d}=\prod_{n=0}^\infty\left(\sum_{j=0}^\infty a_n^{-dj}\right)=\prod_{n=0}^\infty\frac{1}{1-a_n^{-d}}
\]
which converges to a non-zero number in view of \eqref{super}.
\end{proof}
    We set for $p=1,2,\ldots$
      \begin{equation}
        \label{f-p}
     \mathcal F_{-p}=   \left\{(f_\alpha)_{\alpha\in\ell}\,\,; \sum_{\alpha\in\ell}|f_\alpha|^2 b_\alpha^{-p}<\infty\right\}
      \end{equation}
and
      \begin{equation}
        \label{f-p-1}
        \mathcal F_{p}=   \left\{(f_\alpha)_{\alpha\in\ell}\,\,; \sum_{\alpha\in\ell}|f_\alpha|^2
          (\alpha!)^2b_\alpha^{p}<\infty\right\}.
      \end{equation}
The space of stochastic test functions is
\begin{equation}
        \mathfrak F_{test}(a)=\bigcap_{p=1}^\infty \mathcal F_{p}
        \label{Fa-1}
        \end{equation}
and the space of stochastic distributions is defined by
\begin{equation}
        \mathfrak F_{dist}(a)=\bigcup_{p=1}^\infty \mathcal F_{-p}.
        \label{Fa}
      \end{equation}
      We have
      \[
        \mathfrak F_{test}(a)\subset \mathbf L_2(\mathscr S_{\mathbb R}^\prime,\mathcal C,P)
        \subset \mathfrak F_{dist}(a).
\]
      By the analysis in \cite{GS2_english,MR0435834}, the space $\mathfrak F_{test}(a)$ is nuclear Fr\'echet, and is in duality with $\mathfrak F_{dist}(a)$ via
\begin{equation}
  \label{erty}
\langle f,g\rangle=\sum_{\alpha\ell} \alpha!\overline{g_\alpha}f_\alpha.
\end{equation}
For completeness we prove:
\begin{lem}
  \label{erty1}
\eqref{erty} makes sense for all $f\in\mathfrak F_{test}(a)$ and $g\in\mathfrak F_{dist}(a)$.
\end{lem}

\begin{proof}
  Let $p\in\mathbb N$ be such that $g\in\mathcal  F_{-p}$. By definition of $\mathfrak F_{test}(a)$
  we have $f\in\mathcal F_{p}$ and so by the Cauchy-Schwarz inequality we have:
  \[
    \begin{split}
      |\sum_{\alpha\ell} \overline{g_\alpha}(\alpha!)^2f_\alpha|&=|\sum_{\alpha\ell} \overline{g_\alpha}2^{-p/2}(\alpha!)2^{p/2}f_\alpha|\\
      &\le\left(\sum_{\alpha\in\ell}|g_\alpha|^22^{-np}\right)^{1/2}\left(\sum_{\alpha\in\ell}|f_\alpha|^2(\alpha!)^22^{np}\right)^{1/2}<\infty.
    \end{split}
    \]
    \end{proof}
      
      \begin{thm}  (see \cite{MR3029153,MR3404695})
        With the convolution of coefficient as product, if $f\in\mathcal F_{-p}$ and $g\in\mathcal F_{-q}$ with $p-q\ge d$, then
        $f\star g\in\mathcal F_{-p}$        and
        \begin{equation}
          \|f\star g\|_{-p}\le A(p-q)\|f\|_{-p}\cdot\|g\|_{-q}.
        \end{equation}
        \label{dallas}
        \end{thm}

        \begin{lem}
\label{f2-nucl}
          In the previous notation and hypothesis, the space $\mathfrak F_{dist}(a)$ is nuclear.
\end{lem}

\begin{proof}
This is \cite[Corollary 3.4]{vage1} applied to the sequence $(2^{n})$ and $p=2$.
  \end{proof}

  \section{The complex Brownian motion $B_z$}
  \label{sec-5}
\setcounter{equation}{0}
Let $R>0$ and let $C_R$ be defined by \eqref{C_R}. Let furthermore $|z|\le R$. Since $|y|\le R$ it follows from \eqref{ineq-k} that
\[
|\zeta_n(z)|\le C_R\cdot e^{2Rn}
  \]
  so that
  \begin{equation}
    \left|\int_{[0,z]}\zeta_n(u)du\right|\le R C_R\cdot e^{2Rn}.
  \end{equation}

  \begin{defn}
    \label{2}
  We set $a_n=2^n$, $n=0,1,2,\ldots$, and denote by $\mathbf 2$ the corresponding sequence. We trivially have \eqref{super1} and \eqref{super} holds for $d=1$. We let $\mathfrak F_{dist}(\mathbf 2)$ be the corresponding strong algebra as in
  \eqref{Fa} and $\mathfrak F_{test}(\mathbf 2)$ be the corresponding space of test functions.
\end{defn}
  \begin{thm}
    \label{bounded-thm}
    The map $z\mapsto B_z$ is continuous form $\overline{B}(0,R)$ into $\mathfrak F_{dist}(\mathbf 2)$, and
    there  is a $p$ such that $(B_z)\in\mathcal F_{-p}$ for $z\in \overline{B}(0,R)$ and $B_z$ is bounded in the norm of $\mathcal F_{-p}$ in
    $\overline{B}(0,R)$.
  \end{thm}

  \begin{proof}
    The space $\overline{B}(0,R)$ is a metric space, and so it is enough to work with sequences, which will allow us to take advantage of
    Theorem \ref{gelfand-1}.  Since $\zeta_n$ is an entire function we have
    \[
\int_{[0,z_1]}\zeta_n(u)du+\int_{[z_1,z_2]}\zeta_n(u)du=\int_{[0,z_2]}\zeta_n(u)du,
\]
and so
    \[
      \begin{split}
        |\int_{[0,z_1]}\zeta_n(u)du-\int_{[0,z_2]}\zeta_n(u)du|&=|\int_{[z_2,z_1]}\zeta_n(u)du|\\
        &\le |z_1-z_2|\max_{z\in \overline{B}(0,R)}|\zeta_n(u)|\\
        &\le |z_1-z_2| RC_R\cdot e^{2Rn}
\end{split}
\]
and hence
      \[
        \begin{split}
\|B_{z_1}-B_{z_2}\|_{-p}^2&\le R^2(C_R)^2 |z_1-z_2|^2 \left(\sum_{n=0}^\infty \left(e^{2Rn} \right)^2(2^n)^{-p}\right)<\infty
          \end{split}
        \]
        as soon as $p$ is such that
        \begin{equation}
          \label{bound-p}
\frac{e^{4R}}{2^p}<1.
          \end{equation}

        Setting $z_2=0$ and $z_1=z$ in the above it follows that $(B_z)$ belong to $\mathcal F_{-p}$, with $p$ satisfying \eqref{bound-p}.
        Using sequences and Theorem \ref{gelfand-1} we conclude that the function $z\mapsto B_z$ is continuous from $\overline{B}(0,R)$ into
        $\mathfrak F_{dist}(\mathbf 2)$ and bounded in norm there.\smallskip
    
    \end{proof}

    \section{The complex white noise $N_z$}
    \label{white-z}
    \setcounter{equation}{0}

    We define
    \begin{equation}
      N_z(\w)=\sum_{n=0}^\infty \zeta_n(z) Z_n(\w)
    \end{equation}
    where the $Z_n$ have been defined in \eqref{zn}.
    \begin{thm}
      \label{analyticity-thm}
            Let $R>0$ and let $|z|\le R$. Then, $N_z\in\mathcal F(a)$. Mroe precisely, there is a $p$ independent of $|z|\le R$ such that
            $N_z\in\mathcal F_{-p}$. Moreover,  in the topology of
            $\mathcal F(a)$
            \begin{equation}
              \partial B_z=N_z.
              \end{equation}
            \end{thm}
            \begin{proof}
              Using \eqref{ineq-k-2} to go from the last-to-one line to the last one, we have
              \[
                \begin{split}
                  \left|\frac{\int_{[0,z]}\zeta_n(u)du-\int_{[0,z_0]}\zeta_n(u)du}{z-z_0}-\zeta_n(z_0)\right|
                  &=\left|\frac{\int_{[z_0,z]}(\zeta_n(u)-\zeta_n(z_0))du}{z-z_0}\right|\\
                  &=\left|\frac{\int_0^1(\zeta_n(z_0+t(z-z_0))-\zeta_n(z_0))(z-z_0)dt}{z-z_0}\right|\\
                  &\le \max_{t\in[0,1]}|\zeta_n(z_0+t(z-z_0))-\zeta_n(z_0)|\\
                  &= \max_{t\in[0,1]}|\int_{[z_0,z_0+t(z-z_0)]}\zeta_n^\prime(w)dw|\\
                                    &= \max_{t\in[0,1]}|\int_0^1\zeta_n^\prime(z_0+st(z-z_0))t(z-z_0)ds|\\
                  &\le |z-z_0|\cdot \max_{u\in[z_0,z]}|\zeta_n^\prime(u)|\\
                  &\le |z-z_0|\cdot C_R                  \left(1+\sqrt{2n}\right)e^{2n|R|}.
                  \end{split}
\]
Thus, with $p$ satisfying
\begin{equation}
  \label{bound-p-2}
  \frac{2ne^{4nR}}{2^p}<1
\end{equation}
we have
              \[
                \begin{split}
                  \left\|\frac{B_{z}-B_{z_0}}{z-z_0}-N_{z_0}\right\|_p^2&=\sum_{n=0}^\infty
                  \left|\frac{\int_{[0,z]}\zeta_n(u)du-\int_{[0,z_0]}\zeta_n(u)du}{z-z_0}-\zeta_n(z_0)\right|^2\\
                  &\le C_R^2 |z-z_0|^2\left(\sum_{n=0}^\infty (\left(1+\sqrt{2n}\right)e^{2n|R|})2^{-pn}\right)\\
                      &                                  \rightarrow 0\quad{\rm as}\,\,\,|z-z_0|\rightarrow 0
                \end{split}
              \]
              So the convergence occurs in particular via sequences and there is convergence in $\mathcal F_{-p}$ with $p$ satisfying
              \eqref{ineq-k-2}.
            \end{proof}

The following result is implicit in the proof of the previous theorem, but we detail the arguments in the proof of the following result:
            
\begin{thm}
  \label{nz-conti}
The map $z\mapsto N_z$ is continuous from $\mathbb C$ into $\mathfrak F_{dist}(\mathbf 2)$.
\end{thm}

\begin{proof}
  Let $n\in\mathbb N_0$, and $z,w\in\overline{B}(0,R)$. Using \eqref{ineq-k-3} and the definition of $C_R$ (see \eqref{K_z}--\eqref{C_R}) we have
  \[
    \begin{split}
      |\zeta_n(z)-\zeta_n(w)|&=|\int_{[w,z]}\zeta^\prime (v)dv|\\
      &\le|z-w|\max_{v\in[z,w]}|\zeta^\prime(v)|\\
      &\le |z-w|C_R(1+\sqrt{2n})e^{2nR}.
      \end{split}
    \]
    Thus
    \[
\|\zeta_n(z)-\zeta_n(w)|^2_{-p}\le |z-w|^2\cdot C_R^2\left(\sum_{n=0}^\infty\frac{(1+\sqrt{2n})^2e^{4nR}}{2^{np}}\right)
\]
where the infinite sum will converge for any $p$ satisfying \eqref{bound-p-2} (and in particular \eqref{bound-p}).
It follows that $z\mapsto N_z$ is continuous from the convex compact set $\overline{B}(0,R)$ into $\mathcal F_{-p}$ for such values of $p$ and hence, by Theorem
\ref{gelfand-1} from $\overline{B}(0,R)$ into $\mathcal F_{-p}$ into $\mathfrak F_{dist}(\mathbf 2)$.
\end{proof}

\begin{rem}
In the preceding argument sequences are not used to prove the continuity.
\end{rem}

The proof of Theorem \ref{analyticity-thm} implies in fact much more:

\begin{thm}
  \label{6.4}
The map $z\mapsto B_z$ is analytic.
\end{thm}

\begin{proof}
  Let $R>0$. By Theorem \ref{bounded-thm} there exists a $p\in\mathbb N$ such that $B_z\in\mathcal F_{-p}$ and is uniformly bounded there in norm for
  $|z|\le R$. The latter property insures that we can use \cite[Proposition A.3 p. 462]{MR2798103}. To prove analyticity is equivalent to prove that, for every
  $f\in\mathcal F_p$ the function $z\mapsto \langle B_z,f\rangle$ is analytic, where the brackets denote the duality \eqref{erty}.
  But
  \[
\langle B_z,f\rangle=\sum_{n=0}^\infty\overline{f_{\alpha^{(n)}}}\int_{[0,z]}\zeta_n(u)du
    \]
    where $\alpha^{(n)}$ is the element of $\ell$ with all entries equal to $0$, besides the $n$-th one, equal to $1$. By Lemma \ref{erty1} this series is a uniformly convergent
    series of analytic functions, and hence is analytic.\smallskip

    At this space we have shown that $B_z$ is analytic from $|z|<R$ to $\mathcal F_{-p}$, in the topology of $\mathcal F_{-p}$. By definition of the inductive
    topology, this implies analyticity from $B_z$ is analytic from $|z|<R$ into $\mathfrak F_{dist}(a)$.
  \end{proof}
              \section{Stochastic integral}
              \setcounter{equation}{0}
              \label{sec-7}
              Given a function $f(z)$ continuous from $\overline{B}(0,R)$ into $\mathfrak F_{dist}(\mathbf 2)$ we wish to define the stochastic integral $\int_C f(z)\star N_z dz$
              where $C$ is a smooth path inside $\overline{B}(0,R)$. Taking a continuously differentiable parametrization $\gamma(t),\, t\in[a,b]$ of $C$ we wish to define
            \begin{equation}
          \label{integr}
          \int_C f(z)\star N_z dz=\int_a^bf(\gamma(t))\star N_{\gamma(t)}\gamma^\prime(t)dt.
        \end{equation}
        Of course the integral, if it makes sense at all, will depend on the path when $f$ is merely assumed continuous.\\

In the following theorem  we prove that the function
            \begin{equation}
              \label{f-n-g}
  t\mapsto f(\gamma(t))\star N_{\gamma(t)}\gamma^\prime(t),\quad t\in[a,b],
\end{equation}
is continuous from $[a,b]$ into one of the spaces $\mathcal F_{-p}$. This allows to define \eqref{integr}
as a Riemann integral of a continuous function in a Hilbert space.

\begin{thm}
  \label{7.1}
Let $f(z)$ be continuous from $\overline{B}(0,R)$ into $\mathfrak F_{dist}(\mathbf 2)$.
The integral \eqref{integr} exists the integral of a Hilbert space valued integral of a continuous function.
        \end{thm}

        \begin{proof}
          We proceed in a number of steps.\smallskip

          STEP 1: {\sl The Wick product is jointly continuous in the two variables in $\mathfrak F_{dist}(\mathbf 2)$.}\smallskip

          Indeed,  the nuclearity of $\mathfrak F_{dist}(\mathbf 2)$ (see Lemma \ref{f2-nucl})
          implies that a set is bounded there if and only if it is bounded in one of the $\mathcal F_{-p}$
          (see Theorem \ref{gelfand-1}). Thus, by \cite[Theorem 3.3]{MR3404695}, the Wick product is jointly continuous.\smallskip

          STEP 2: {\sl The function \eqref{f-n-g} is continuous from $[a,b]$ into $\mathfrak F_{dist}(\mathbf 2)$.}\smallskip

By Theorem \ref{nz-conti}, $z\mapsto N_z$ is continuous from $\mathbb C$ into $\mathfrak F_{dist}(\mathbf 2)$ and so by the previous step,
the function $z\mapsto f(z)\star N_z$ is continuous in $\mathfrak F_{dist}(\mathbf 2)$. Thus the function \eqref{f-n-g} is continuous from $[a,b]$ into
$\mathfrak F_{dist}(\mathbf 2)$ since $t\mapsto\gamma(t)$ and $t\mapsto \gamma^\prime(t)$ are continuous by assumption.\smallskip

          STEP3: {\sl The integral \eqref{integr} is a strong Hilbert space valued integral.}\smallskip

          Since $[0,1]$ is compact, and since the continuous image of a compact is compact, the previous step implies that the range of the function
\eqref{f-n-g} is compact in $\mathfrak F_{dist}(\mathbf 2)$. Hence, by the third item in Theorem \ref{gelfand-1}
          it is compact in one of the spaces $\mathcal F_{-p}$, say $\mathcal F_{p_0}$, in the corresponding topology. The integral
          \[
\int_a^b f(\gamma(t))\star N_{\gamma(t)}\gamma^\prime(t)dt
            \]
            exists as a strongly converging Riemann integral in $\mathcal F_{p_0}$. It is the limit with respect to a net. Since the limit exists, we can now take a sequence of Riemann
            sums converging to the integral, and we get convergence in the strong topology of $\mathfrak F_{dist}(\mathbf 2)$ by item four of Theorem \ref{gelfand-1}.\smallskip

            STEP 4: {\sl In the previous step, the limit in $\mathfrak F_{dist}(\mathbf 2)$ will not depend on the chosen sequence.}\smallskip

            Let the integral be computed in the spaces $\mathcal F_{-p_1}$ and $\mathcal F_{-p_2}$ respectively, with limits $\ell_1\in\mathcal F_{-p_1}$ and $\ell_2
\in\mathcal F_{-p_2}$ respectively. We assume $p_1<p_2$ and thus $\mathcal F_{-p_1}\subset \mathcal F_{-p_2}$, with the converse inequality on their norms.
Denonog by $g$ the function \eqref{f-n-g} we have (with obvious notation for the Riemann sums)
\[
\|\left(\sum g(t_j)(s_{j+1}-s_j)\right)-\ell_1\|_{-p_2}\le \|\left(\sum g(t_j)(s_{j+1}-s_j)\right)-\ell_1\|_{-p_1}.
    \]
    Since the net defined by the right hand silde converges so does the one on the left, and so the limit of the integral in $\mathcal F_{-p_2}$ is also $\ell_1$ and so
    $\ell_1=\ell_2$ by uniqueness of the limit and since $\mathcal F_{-p_1}\subset\mathcal F_{-p_2}$.
  \end{proof}

  \section{Real versus imaginary time}
\setcounter{equation}{0}\label{section-8}
\begin{thm}
  Let $T\not=0$. Then, $B_{iT}$ does not belong to $\mathbf L_2(\mathcal S_{\mathbb R}^\prime,\mathcal C,P)$.
\end{thm}

\begin{proof}
  Using the dominated convergence theorem, and since the Hermite functions are real, $\overline{\zeta_n(z)}=\zeta_n(\overline{z})$
  and satisfy \eqref{ineq-k}, we have
\begin{equation}
  \begin{split}
    \iint_{[0,z]\times [0,z]}\left(    \sum_{n=0}^\infty \e^n\zeta_n(u)\overline{\zeta_n(v)}\right)dudv&=\\
    &\hspace{-2cm}=
    \sum_{n=0}^\infty \e^n\iint_{[0,1]\times [0,1]}\zeta_n(zt){\zeta_n(\overline{z}t^\prime)}z\overline{z}dtdt^\prime\\
    &\hspace{-2cm}=    \sum_{n=0}^\infty \e^n\left|\int_{[0,z]}\zeta_n(v)dv\right|^2\\
    &\hspace{-2cm}<\infty
  \end{split}
  \label{1}
\end{equation}
for $\e$ satisfying
\[
  |\e|\cdot e^{{2}|z|}<1.
  \]
On the other hand  Mehler's formula gives us
\[
   \sum_{n=0}^\infty \e^n\zeta_n(zt){\zeta_n(\overline{z}t^\prime)}
 =\frac{1}{\sqrt{\pi}}\frac{1}{\sqrt{1-\e^2}}e^{-\frac{(1+\e^2)(z^2t^2+\overline{z}^2
      (t^\prime)^2)-4|z|^2tt^\prime\e}{2(1-\e^2)}}.
\]
Take $z=iT$ with $T\in\mathbb T$. Since $t,t^\prime\in[0,1]$ the above equality becomes 
\[
  \begin{split}
   \sum_{n=0}^\infty \e^n\zeta_n(iTt){\zeta_n(-iTt^\prime)}
   & =\frac{1}{\sqrt{\pi}}\frac{1}{\sqrt{1-\e^2}}e^{-\frac{-(1+\e^2)(T^2(t^2+    (t^\prime)^2)-4|T|^2tt^\prime\e}{2(1-\e^2)}}\\
&   =\frac{1}{\sqrt{\pi}}\frac{1}{\sqrt{1-\e^2}}e^{\frac{(1+\e^2)(T^2(t^2+    (t^\prime)^2)+4|T|^2tt^\prime\e}{2(1-\e^2)}}\\
   &\ge \frac{1}{\sqrt{\pi}}\frac{1}{\sqrt{1-\e^2}}e^{\frac{(1+\e^2)(T^2(t^2+    (t^\prime)^2))}{2(1-\e^2)}}\\
   \end{split}
\]
and so
\[
  |z|^2\sum_{n=0}^\infty\e^n\left|\int_0^1\zeta_n(zt)dt\right|^2\ge\frac{1}{\sqrt{\pi}}\frac{1}{\sqrt{1-\e^2}}\left(\int_0^1 e^{\frac{(1+\e^2)T^2t^2}{2(1-\e^2)}}dt
    \right)^2
  \]
  The monotone convergence theorem  shows that
  \begin{equation}
    \label{ineq345}
    \lim_{\e\uparrow 1}\sum_{n=0}^\infty\e^n\left|\int_0^1h(zt)dt\right|^2=\sum_{n=0}^\infty \left|\int_0^1\zeta_n(zt)dt\right|^2
\end{equation}
and it follows from \eqref{ineq345} that $B_z$ is not in
the white noise space for $z=iT$, $T\in\mathbb R\setminus\left\{0\right\}$.
\end{proof}

        \section{A more general family of processes}
        \setcounter{equation}{0}
        We extend the previous analysis to Gaussian processes $(X_t)_{t\ge 0}$  indexed by $[0,\infty)$ and with covariance
        \begin{equation}
          \label{cov-m}
E(X_tX_s)=\int_0^{t\wedge s}m(u)du,
\end{equation}
where $m(u)$ has an analytic squareroot  in $\mathbb C$ (and is in particular analytic there), and positive on $[0,\infty)$. Examples of such functions $m$ include $e^{\pm u^{2n}}$ (with $R=\infty$) and $m(u)=(1+u^2)^{2n}$.

        \begin{lem}
          \begin{equation}
X_t(\w)=\sum_{n=0}^\infty \left(\int_0^t\sqrt{m(u)}\zeta_n(u)du\right)Z_n(\w).
\end{equation}
is centered Gaussian, with covariance function \eqref{cov-m}.
          \end{lem}
          \begin{proof}
            We have
            \begin{equation}
1_{[0,t]}(u)\sqrt{m(u)}\zeta_n(u)=\sum_{n=0}^\infty \left(\int_0^t\sqrt{m(v)}\zeta_n(v)dv\right)\zeta_n(u).
              \end{equation}
            where the convergence is in $\mathbf L_2(\mathbb R,du)$.
          \end{proof}

          We define
          \begin{equation}
X_z=\sum_{n=0}^\infty \left(\int_{[0,z]}\sqrt{m(u)}\zeta_n(u)du\right)Z_n(\w),
\end{equation}
where $\sqrt{m}$ denotes the analytic squareroot of $m$ in $B(0,R)$ which coincide with $\sqrt{m(u)}$ on $(0,R)$
The coefficients of the process $(X_z)$ have the same bound as $B_z$ since $m(z)$ is bounde don compacts.

\begin{thm}
  \label{9.2}
  The process $(X_z)_{z\in\mathbb C}$ is differentiable (i.e. analytic) in $\mathbb C$, with derivative
  \begin{equation}
    \frac{\partial X_z}{\partial z}(z)=\sqrt{m(z)}N_z(\w)
  \end{equation}
  converging uniformly on compact subsets of the complex plane.
  \end{thm}

  \begin{proof}
    The argument for the Brownian motion and its derivative still hold here since $\sqrt{m(z)}$ is bounded in modulus on compact sets. Furthermore, the derivative is now
\[
    \sqrt{m(z)}\sum_{n=0}^\infty \zeta_n(z)Z_n(\w)=\sqrt{m(z)}N_z(\w).
\]
    \end{proof}
    \section{An extension of It\^o’s formula to the case of complex time: The quadratic case}
    \setcounter{equation}{0}

    We recall that
    \begin{equation}
      \label{equ-10-1}
      Z_n\star Z_m=\begin{cases}Z_n Z_m,\quad n\not=m\\
        Z_n^2-1,\quad n=m
        \end{cases}
      \end{equation}
      Using these identiy and for {\bf real}\, $t$ we have a formal proof of
    \begin{equation}
      B_t^2=\int_{0}^t 2B_u\star N_u dt +t
    \end{equation}
    as follows:
    \[
      \begin{split}
        B_t^2-\int_{0}^t 2B_t\star N_t dt&=\\
        &\hspace{-20mm}=\sum_{n,m=0}^\infty\left(\int_0^t\zeta_n(u)du\right)\left(\int_0^t\zeta_m(u)du\right) Z_n Z_m-\\
        &\hspace{-12.mm}-2\int_0^t\left(\sum_{n,m=0}^\infty\left(\int_0^u\zeta_n(v)dv\right)\zeta_m(u)du\right)Z_n\star Z_m\\
                &\hspace{-20mm}=\sum_{n,m=0}^\infty\left(\int_0^t\zeta_n(u)du\right)\left(\int_0^t\zeta_m(u)du\right) Z_n Z_m-\\
                &\hspace{-12.mm}-2\sum_{n,m=0}^\infty\left(\int_0^t\left(\int_0^u\zeta_n(v)dv\right)\zeta_m(u)du\right)Z_n\star Z_m\\
                &\hspace{-20mm}=\sum_{n,m=0}^\infty\left(\int_0^t\zeta_n(u)du\right)\left(\int_0^t\zeta_m(u)du\right) (Z_n Z_m-\delta_{n,m}+\delta_{n,m})-\\
                &\hspace{-12.mm}-2\sum_{n,m=0}^\infty\left(\int_0^t\left(\int_0^u\zeta_n(v)dv\right)\zeta_m(u)du\right)Z_n\star Z_m\\
                &\hspace{-20mm}=\sum_{n,m=0}^\infty\left(\int_0^t\zeta_n(u)du\right)\left(\int_0^t\zeta_m(u)du\right) (Z_n \star Z_m+\delta_{n,m})-\\
                &\hspace{-12.mm}-2\sum_{n,m=0}^\infty\left(\int_0^t\left(\int_0^u\zeta_n(v)dv\right)\zeta_m(u)du\right)Z_n\star Z_m\\
                &\hspace{-20mm}=\underbrace{\left(\sum_{n,m=0}^\infty\left(\int_0^t\zeta_n(u)du\right)\left(\int_0^t\zeta_n(u)du\right)\right)}_{t\,\,
                  {\rm by \,\, Parseval}}
                +\sum_{n,m=0}^\infty A_{n,m} Z_n\star Z_m\\
                &\hspace{-20mm}=t+\sum_{n,m=0}^\infty A_{n,m} Z_n\star Z_m
      \end{split}
    \]
    where
    \[
      A_{n,m}=\left(\int_0^t\zeta_n(v)dv\right)\left(\int_0^t\zeta_m(u)du\right)-      \int_0^t\left(\int_0^u\zeta_n(v)dv\right)\zeta_m(u)du.
      \]
    But
    \[
      \int_0^t\left(\int_0^u\zeta_n(v)dv\right)\zeta_m(u)du=\left(\int_0^t\zeta_n(v)dv\right)\left(\int_0^t\zeta_m(u)du\right)-
      \int_0^t\zeta_n(u)\int_0^u\zeta_m(v)dv,
    \]
    that is
 \[
      \int_0^t\left(\int_0^u\zeta_n(v)dv\right)\zeta_m(u)du+\int_0^t\zeta_n(u)\int_0^u\zeta_m(v)dv=\left(\int_0^t\zeta_n(v)dv\right)\left(\int_0^t\zeta_m(u)du\right)      
            \]
            and so $A_{n.m}+A_{m,n}=0$, and the result follows.\\

            To adapt this computation to the present setting we will replace $B_t$ by
            \[
B_{z,\e}=\sum_{n=0}^\infty \e^{n/2}\left(\int_{[0,z]}\zeta_n(u)du\right)Z_n
\]
where $\e$ is such that $|\e|<1$.
By the bound \eqref{ineq-k} we have that $B_{\e,z}$ belongs to the white noise space. We will
            will require the counterpart of Parseval equality. 
            Mehler's formula and 
            We have \eqref{mehler678}, i.e.
            \begin{equation*}
\iint_{[0,z]\times[0,z]}  \left(   \sum_{n=0}^\infty\e^n\zeta_n(u)\zeta_n(v)  \right)dudv= \frac{1}{\sqrt{2\pi}}\frac{1}{\sqrt{1-\e^2}}\iint_{[0,z]\times[0,w]}e^{-\frac{(1+\e^2)(u^2+\overline{v}^2)-4u\overline{v}\e}
        {2(1-\e^2)} }     dudv
    \end{equation*}
    which is the counterpart of
    \[
      \sum_{n=0}^\infty\left(\int_0^t\zeta_n(u)du\right)\left(\int_0^t\zeta_n(u)du\right)=t.
    \]
    It\^o's formula is now
    \[
B_{z,\e}^2=\int_{[0,z]} 2B_{u,\e}\star N_u du+  \frac{1}{\sqrt{2\pi}}\frac{1}{\sqrt{1-\e^2}}\iint_{[0,z]\times[0,w]}e^{-\frac{(1+\e^2)(u^2+\overline{v}^2)-4u\overline{v}\e}
        {2(1-\e^2)} }     dudv.
      \]

\bibliographystyle{plain}

\begin{thebibliography}{10}

\bibitem{MR3231624}
D.~Alpay, P.~Jorgensen, and G.~Salomon.
\newblock On free stochastic processes and their derivatives.
\newblock {\em Stochastic Process. Appl.}, 124(10):3392--3411, 2014.

\bibitem{MR2610579}
D.~Alpay and D.~Levanony.
\newblock Linear stochastic systems: a white noise approach.
\newblock {\em Acta Appl. Math.}, 110(2):545--572, 2010.

\bibitem{vage1}
D.~Alpay and G.~Salomon.
\newblock New topological {$\Bbb C$}-algebras with applications in linear
  systems theory.
\newblock {\em Infin. Dimens. Anal. Quantum Probab. Relat. Top.},
  15(2):1250011, 30, 2012.

\bibitem{MR3029153}
D.~Alpay and G.~Salomon.
\newblock Topological convolution algebras.
\newblock {\em J. Funct. Anal.}, 264(9):2224--2244, 2013.

\bibitem{MR3404695}
D.~Alpay and G.~Salomon.
\newblock On algebras which are inductive limits of {B}anach spaces.
\newblock {\em Integral Equations Operator Theory}, 83(2):211--229, 2015.

\bibitem{MR2798103}
W.~Arendt, C.J.K. Batty, M.~Hieber, and F.~Neubrander.
\newblock {\em Vector-valued {L}aplace transforms and {C}auchy problems},
  volume~96 of {\em Monographs in Mathematics}.
\newblock Birkh\"{a}user/Springer Basel AG, Basel, second edition, 2011.

\bibitem{MR3456934}
X.~Bardina and C.~Rovira.
\newblock Approximations of a complex {B}rownian motion by processes
  constructed from a {L}\'{e}vy process.
\newblock {\em Mediterr. J. Math.}, 13(1):469--482, 2016.

\bibitem{MR0228020}
W.~Feller.
\newblock {\em An introduction to probability theory and its applications.
  {V}ol. {I}}.
\newblock Third edition. John Wiley \& Sons, Inc., New York-London-Sydney,
  1968.

\bibitem{GS2_english}
I.M. Gelfand and G.E. Shilov.
\newblock {\em Generalized functions. Volume 2}.
\newblock Academic Press, 1968.

\bibitem{MR0435834}
I.M. Gelfand and N.Ya. Vilenkin.
\newblock {\em Generalized functions. {V}ol. 4}.
\newblock Academic Press [Harcourt Brace Jovanovich Publishers], New York, 1964
  [1977].
\newblock Applications of harmonic analysis, Translated from the Russian by
  Amiel Feinstein.

\bibitem{MR4666291}
E.~Gwynne, J.~Pfeffer, and M.~Park.
\newblock Loewner evolution driven by complex {B}rownian motion.
\newblock {\em Ann. Probab.}, 51(6):2086--2130, 2023.

\bibitem{Hida_BM}
T.~Hida.
\newblock {\em Brownian motion}, volume~11 of {\em Applications of
  Mathematics}.
\newblock Springer-Verlag, New York, 1980.
\newblock Translated from the Japanese by the author and T. P. Speed.

\bibitem{MR562914}
T.~Hida.
\newblock {\em Brownian motion}, volume~11 of {\em Applications of
  Mathematics}.
\newblock Springer-Verlag, New York-Berlin, 1980.
\newblock Translated from the Japanese by the author and T. P. Speed.

\bibitem{Hi90}
T.~Hida.
\newblock Functionals of {B}rownian motion.
\newblock In {\em Lectures in applied mathematics and informatics}, pages
  286--329. Manchester Univ. Press, Manchester, 1990.

\bibitem{MR1502747}
E.~Hille.
\newblock A class of reciprocal functions.
\newblock {\em Ann. of Math. (2)}, 27(4):427--464, 1926.

\bibitem{new_sde}
H.~Holden, B.~{\O}ksendal, J.~Ub{\o}e, and T.~Zhang.
\newblock {\em Stochastic partial differential equations}.
\newblock Universitext. Springer, New York, second edition, 2010.
\newblock A modeling, white noise functional approach.

\bibitem{MR0205028}
J.~Horv{\'a}th.
\newblock {\em Topological vector spaces and distributions. {V}ol. {I}}.
\newblock Addison-Wesley Publishing Co., Reading, Mass.-London-Don Mills, Ont.,
  1966.

\bibitem{MR1851117}
Zh. Huang and J.~Yan.
\newblock {\em Introduction to infinite dimensional stochastic analysis},
  volume 502 of {\em Mathematics and its Applications}.
\newblock Kluwer Academic Publishers, Dordrecht, chinese edition, 2000.

\bibitem{MR1407325}
Y.G. Kondratiev, P.~Leukert, and L.~Streit.
\newblock Wick calculus in {G}aussian analysis.
\newblock {\em Acta Appl. Math.}, 44(3):269--294, 1996.

\bibitem{MR1387829}
H.-H. Kuo.
\newblock {\em White noise distribution theory}.
\newblock Probability and Stochastics Series. CRC Press, Boca Raton, FL, 1996.

\bibitem{lifs}
M.A. Lifshits.
\newblock {\em Gaussian random functions}, volume 322 of {\em Mathematics and
  its {A}pplications}.
\newblock Kluwer {A}cademic {P}ublisher, 1995.

\bibitem{MR521027}
A.~G. Miamee and H.~Salehi.
\newblock On the factorization of a nonnegative operator valued function.
\newblock In {\em Probability theory on vector spaces (Proc. Conf.,
  Trzebieszowice, 1977)}, volume 656 of {\em Lecture Notes in Math.}, pages
  129--137. Springer, Berlin, 1978.

\bibitem{MR161189}
E.~Nelson.
\newblock Feynman integrals and the {S}chr\"{o}dinger equation.
\newblock {\em J. Mathematical Phys.}, 5:332--343, 1964.

\bibitem{MR214150}
E.~Nelson.
\newblock {\em Dynamical theories of {B}rownian motion}.
\newblock Princeton University Press, Princeton, NJ, 1967.

\bibitem{MR1808423}
P.~Rusev.
\newblock An inequality for {H}ermite's polynomials in the complex plane.
\newblock {\em C. R. Acad. Bulgare Sci.}, 53(10):13--16, 2000.

\bibitem{MR3882989}
L.~Sang, G.~Shen, and Q.~Wang.
\newblock Weak convergence of the complex fractional {B}rownian motion.
\newblock {\em Adv. Difference Equ.}, pages Paper No. 444, 18, 2018.

\bibitem{sansone}
G.~Sansone.
\newblock {\em Orthogonal functions}.
\newblock Dover {P}ublications, {I}nc., {N}ew--{Y}ork, 1991.
\newblock Revised {E}nglish {E}dition.

\bibitem{szego}
G.~Szeg{\"{o}}.
\newblock {\em Orthogonal polynomials}, volume~23.
\newblock {Amer.} Math. Soc., {Rhodes} Island, 1975.

\bibitem{thangavelu1993lectures}
S.~Thangavelu.
\newblock {\em Lectures on Hermite and Laguerre expansions}, volume~42.
\newblock Princeton University Press, 1993.

\bibitem{MR37:726}
F.~Tr{\`e}ves.
\newblock {\em Topological vector spaces, distributions and kernels}.
\newblock Academic Press, New York, 1967.

\bibitem{4839872}
user12030145 (https://math.stackexchange.com/users/1114056/user12030145).
\newblock A uniform (non-asymptotic) upper bound for {H}ermite polynomials in
  the complex plane.
\newblock
  \url{https://math.stackexchange.com/questions/4839872/a-uniform-non-asymptotic-upper-bound-for-\\
  Hermite-polynomials-in-the-complex-pl}.

\bibitem{vage96b}
G.~V{\r{a}}ge.
\newblock A general existence and uniqueness theorem for {W}ick-{SDE}s in
  {${S}^n_{-1,k}$}.
\newblock {\em Stochastic {S}ochastic {R}ep.}, 58:259--284, 1996.

\bibitem{vage96}
G.~V{\r{a}}ge.
\newblock Hilbert space methods applied to stochastic partial differential
  equations.
\newblock In H.~K{\"o}rezlioglu, B.~{\O}ksendal, and A.S. {\"U}st{\"u}nel,
  editors, {\em Stochastic analysis and related topics}, pages 281--294.
  Birk{\"a}user, {B}oston, 1996.

\bibitem{MR626346}
N.~N. Vakhania.
\newblock {\em Probability distributions on linear spaces}.
\newblock North-Holland Publishing Co., New York, 1981.
\newblock Translated from the Russian by I. I. Kotlarski, North-Holland Series
  in Probability and Applied Mathematics.

\bibitem{MR1041203}
S.~J.~L. van Eijndhoven and J.~L.~H. Meyers.
\newblock New orthogonality relations for the {H}ermite polynomials and related
  {H}ilbert spaces.
\newblock {\em J. Math. Anal. Appl.}, 146(1):89--98, 1990.

\end{thebibliography}
\def\cprime{$'$} \def\cprime{$'$} \def\cprime{$'$}
  \def\lfhook#1{\setbox0=\hbox{#1}{\ooalign{\hidewidth
  \lower1.5ex\hbox{'}\hidewidth\crcr\unhbox0}}} \def\cprime{$'$}
  \def\cprime{$'$} \def\cprime{$'$} \def\cprime{$'$} \def\cprime{$'$}
  \def\cprime{$'$}

\end{document}